\documentstyle[12pt]{amsart}
\theoremstyle{plain}
\newtheorem{Thm}{ }[section]

\title{Petri map for rank two bundles with canonical
determinant}
\author{montserrat teixidor i bigas}
\address{Mathematics Department, Tufts University, Medford MA
02155}
\email{montserrat.teixidoribigas@@tufts.edu}

\begin{document}
 \begin{abstract} We prove the Bertram-Feinberg-Mukai conjecture
 for a generic curve $C$ of genus $g$ and a semistable vector bundle $E$ of
rank two and determinant $K$ on $C$, namely we prove the injectivity
of the Petri-canonical map $S^2(H^0(E))\rightarrow H^0(S^2(E))$.

\end{abstract}

\maketitle
\begin{section}{Introduction}

Consider a generic curve of genus $g$. Classical Brill-Noether
Theory deals with the existence of line bundles of a fixed degree
$d$ that possess a preassigned number of sections $k$ traditionally
denoted by $r+1$. The loci of such line bundles $L$ is a subvariety
of the Picard scheme $Pic^d(C)$ and is usually denoted by $W^r_d$.

A key role in the study of classical Brill-Noether Theory is
played by the Petri map: Fix a line bundle $L$ such that
$h^0(C,L)=r+1$. Consider the map
$$H^0(C,L)\otimes H^0(C,\omega \otimes L^{-1})\rightarrow
H^0(C,\omega).$$ The orthogonal to the image of this map is
identified with the tangent space to $W^r_d$ at the point $L$. For
generic $C$ and every $L$, this map is injective. This fact alone
proves that $W^r_d$ is empty when its expected dimension is negative
and of expected dimension when non-empty. Moreover, its singular
locus is $W^{r+1}_d$.

The most straightforward generalization of Brill-Noether Theory to
higher rank is to consider the loci inside the moduli space of
stable vector bundles of rank $r$ and degree $d$ that have at least
a given number $k$ of sections.
 The
equivalent of the Petri map is then
$$H^0(C,E)\otimes H^0(C,\omega \otimes E^*)\rightarrow H^0(C,
\omega \otimes E\otimes E^*).$$ The vector space on the right hand
side is identified to the cotangent space to the moduli space of
vector bundles of fixed rank and degree (but variable determinant).
As in the line bundle case, the orthogonal to the image of this map
is the tangent space at $E$ to the locus of vector bundles with $k$
sections. These loci, though, don't enjoy most of the good
properties of classical Brill-Noether theory, in particular the
Petri map is not injective for many stable vector bundles even on a
generic curve (see \cite{duke}).

 On the other hand if one restricts to the case of rank two
 with fixed determinant the canonical line bundle $\omega$ ,
 one obtains the moduli space $U(2,\omega )$. Bertram, Feinberg
 and Mukai (cf.\cite{BF}, \cite{M1}) introduced the
loci $B^k_{2,\omega}$ of (semi)stable vector bundles $E$ of rank $2$
and determinant $\omega$ that have at least $k$ independent
sections.
 They conjectured, based on some evidence for small genus, that
these loci would be well behaved. The object of this paper is to
prove the  part of this conjecture about non-existence and
singularity. The question of existence was partially answered in our
previous work \cite{canonic} as well as in \cite{P}.

We want to study the equivalent of the Petri map in this situation.
When the rank is two and the determinant $\wedge ^2E\cong \omega$,
there is a natural identification of $\omega \otimes E^* \cong E$.
Also, the tangent space to $U(2,\omega)$ inside the moduli space of
stable vector bundles of rank two and degree $2g-2$ can be
identified to the dual of $H^0(C,S^2(E))$ where $S^2(E)$ denotes the
symmetric power of $E$ and is interpreted as a quotient of $E\otimes
E$. As $\wedge ^2H^0(E)$ maps to zero when passing to the quotient,
we can factor the map

 $$H^0(C,E)\otimes H^0(C,E)\rightarrow H^0(C,S^2(E))$$
  by moding out the original space $H^0(C,E)\otimes H^0(C,E)$ by  $\wedge
  ^2H^0(E)$.
    One obtains a map that we
shall call the canonical Petri map or simply the Petri map when
there is no danger of confusion:

$$S^2H^0(C,E)\rightarrow H^0(S^2E).$$
Then the tangent space to $B_{2,\omega}^k$ at $E$ is identified
with the orthogonal to the image of this map. Hence,  its
injectivity shows that the locus is empty when the expected
dimension $\rho =\rho^k_{2,\omega}=3g-3-{k+1\choose 2}<0$, of
dimension exactly $\rho$ when non-empty and singular only along
$B^{k+1}_{2,\omega}$. In particular, when $\rho =-1$, the generic
curve does not have a linear series with so many sections and the
locus of curves where the sections exist is expected to be a
divisor of the moduli space of curves ${\mathcal M}_g$. When
$g=10$ this divisor provided the first counterexample to the slope
conjecture (see \cite{FP} ).  For various values of $k,g$ for
which $\rho =-1$, one can now obtain divisors in ${\mathcal M}_g$
that can provide information on the slope of ${\mathcal M}_g$.

In addition, from work of Mukai \cite{M2}, these loci have been
useful in studying which curves are contained in $K3$ surfaces. For
curves of small genus and for particular values of $k$, their
geometry seems very interesting \cite{OPP} .

Our result is:

\begin{Thm}  {\bf Theorem}
The canonical Petri map is injective for every semistable vector
bundle $E$ of rank two and determinant $\omega$ on a generic curve
of genus $g$ defined over an algebraically closed field of
characteristic different from two.
\end {Thm}
 The proof
is based on  some of the techniques used to show the injectivity of
the classical Petri map in \cite{EH}, although dealing with vector
bundles of rank two is substantially more involved than dealing with
line bundles. The degeneration of the smooth curve to a singular
reducible curve is also different. Instead of a rational curve with
$g$ elliptic tails as in \cite{EH}, we consider  a chain of elliptic
curves (see \cite{W}). This allows us to work in arbitrary
characteristic different from two. Moreover, as pointed out in
\cite{FP} the result fails for at least some of the Eisenbud-Harris
curves.

We will first need  to determine what are the possible limits of
vector bundles with canonical determinant when the curve becomes
reducible. The limit of the vector bundle may become a torsion-free
sheaf for which the determinant is not even defined. We use the
results of Nagaraj-Seshadri (\cite{NS}) as well as the results in
\cite{arbre} to find these limits. Then the theory of limit linear
series for vector bundles (\cite{duke}) can be used to study the
behavior of sections in the kernel.

The author would like to thank the referee for a careful reading of
the manuscript and infinite patience.

\end{section}
\begin{section}{Conditions on semistable vector bundle with canonical
determinant on reducible curves}

Let $C$ be a reducible curve of arithmetic genus $g$ consisting of
two components $C_1,C_2$ with one node obtained by gluing the point
$Q_1\in C_1$ to the point $P_2\in C_2$. Fix a polarization of $C$
(i.e. positive integers $w_1,w_2$ such that $w_1+w_2=1$). Fix a
degree $d$ for rank two bundles and define $\chi =d+2(1-g)$. If
$w_1\chi$ is not an integer, then the moduli space of rank two
bundles on $C$ that are semistable by the given polarization
consists of two components that are characterized by the splitting
of the degrees among $C_1$ and $C_2$ (see \cite{arbre}). Define the
two numbers $\chi _1=\chi (E_{|C_1})$ and $\chi _2=\chi (E_{|C_2})$
where $\chi$ denotes the Euler-Poincar\'e characteristic. Then, the
semistability condition gives
$$(*)w_1\chi \le \chi _1\le w_1\chi +2.$$

We now look at the conditions on the determinant Let $d_1$ be such
that $\chi_1 =d_1+2(1-g)$ is the lowest possible value allowed in
(*). Let $d_2$ be defined as $d_2=d-d_1$. One can then define the
determinant map from the moduli space of torsion free sheaves on $C$
to the Picard variety of line bundles of bidegree $d_1,d_2$ on $C$
as follows (\cite{NS}):

a) If $E$ is a vector bundle of bidegree $d_1,d_2$, take its
determinant.

b) If $E$ is a vector bundle of bidegree $d_1+1,d_2-1$, take
$$((detE_1)(-Q_1), (detE_2)(P_2)).$$

c) If $E$ is a sheaf that is not locally free, then at the node $P$
(see \ref{stable} b) below), $E_P\equiv {\cal O}_P\oplus {\cal
M}_P$. In this case, $w_1\chi <\chi _1<w_1\chi +1$ and $\chi_1 +\chi
_2 =\chi +1$. Hence, the bidegree is $d_1, d_2-1$. The image of $E$
by the determinant map is then
$$(det(E_1), det(E_2(P_2)).$$

These results will be generalized in \ref{detmoltescomp}  to
reducible curves with more components.

 In the rest of the paper we
shall consider a curve of the following type:

\begin{Thm} {\bf Definition}\label{ccgg}
 Assume we have
curves $C_1,..C_M$ each with points $P_i,Q_i$. Glue $Q_i$ to
$P_{i+1}, i=1..M-1$. Assume that every curve is either rational or
elliptic, the points $P_i,Q_i$ are generic on the (elliptic) curve
and the genus of the resulting curve is $g$ (and there are therefore
exactly $g$ elliptic components). We shall call such a curve a chain
of genus $g$
\end{Thm}

Assume that a chain of genus $g$ as above is the central fiber of a
one-dimensional family of curves. If we modify the total space by
base change and blow-ups, the central fiber will still be of the
same kind, may be with a few more rational components.

\begin{Thm} \label{detmoltescomp} {\bf Claim. } Fix  a multidegree $d_1,..d_M$ giving rise to a component
of the moduli space of torsion free sheaves on the total curve $C$
semistable by a certain polarization. Fix another component of the
moduli space of vector bundles of rank two on $C$ corresponding to a
multidegree $d_1+\epsilon_1,..d_M+\epsilon_M$ (hence $\sum \epsilon
_i=0$). Let $E$ be an element of this component with restriction to
$C_i$ being $E_i$ . Then the determinant map takes the form
$$E\rightarrow (det E_i)((\sum _{k<i}\epsilon_k)P_i-(\sum_{k\le
i}\epsilon _k)Q_i)$$
\end{Thm}
\begin{pf}
We use induction on $M$. If $M=2$, this is the result of
Nagaraj-Seshadri about the determinant map. Assume that the formula
is correct for $M-1$ and prove it for $M$. Consider our curve as the
union of two pieces, $C_M$ and the rest. Using the case of two
components, we obtain that on the curve $C_M$ the determinant must
be taken to be $(det E_M)(-\epsilon _M P_M)$. As
$\sum_{i=1...M}\epsilon_i=0$, one finds $\sum
_{k<M}\epsilon_k=-\epsilon_M$ and $\sum_{k\le M}\epsilon _k=0$, this
agrees with the formula above.

On the union of $C_1,..C_{M-1}$, one must take the determinant that
we had in the case of $M-1$ and modify it with
$\epsilon_MQ_{M-1}=(-\sum_{k<M}\epsilon_k)Q_{M-1}$. This
modification does not change anything except on $C_{M-1}$. On
$C_{M-1}$, the restriction of the determinant for the case of $M-1$
components does not have $Q_{M-1}$ in the expression as this is not
a node of $C_{M-1}$. The expression then is $det(E_{M-1})((\sum
_{k<M-1}P_k)$. Modifying this by $\epsilon
_MQ_{M-1}=-(\sum_{k<M}\epsilon_k) Q_{M-1}$ gives the stated result.

\end{pf}

Assume given a family of curves $C_t$ with generic curve
non-singular and special curve $C_0=C$ reducible. Consider a family
of line bundles $L_t$ degenerating to a line bundle $L_0$. Then the
moduli space of vector bundles on $C_t$ with determinant $L_t$
degenerates to the set of torsion free sheaves on $C_0$ with
determinant $L_0$ in the sense above (see \cite{NS}).

In addition to the conditions on the degree of the restriction of
the vector bundle to each component as stated in (*), one finds
conditions on the structure of the restrictions of the vector
bundles themselves. The statement that follows tells us that these
restrictions cannot be too unstable. Its proof is very similar to
the proof of (*) given in \cite{arbre} (see also \cite{G})

\begin{Thm} {\bf Claim}\label{stable} Let $C$ be a curve as in
\ref{ccgg}. Consider torsion-free sheaves of rank two and fixed
Euler-Poincar\'e characteristic $\chi$ and a polarization $\{ w_i\}$
such that $\sum _{i\in I}w_i\chi$ is not an integer for any
$I\subset \{1,\ldots M\} $.

a) The restriction of a semistable sheaf to one of the irreducible
components is either indecomposable or the direct sum of two line
bundles whose degrees differ in at most one unit (or of necessarily
the same degree in the case of the first and last components).

b) If a sheaf is not locally free at one node $P$, then the fiber at
that node is $E_P\cong {\mathcal M}_P\oplus {\mathcal O}_P$ where
${\mathcal M}_P$ denotes the maximal ideal at $P$.
\end{Thm}

\begin{pf} Note that on a rational curve, every vector bundle is a
direct sum of  line bundles while on an elliptic curve a vector
bundle of rank two is either indecomposable or the direct sum of two
line bundles (see \cite{A}). The only statement that needs to be
proved in a) is that the difference between the degrees of these
bundles is either one or zero and that the former is only a
possibility for components with two nodes.

Assume that the restriction of $E$ to a component $C_i$ is the
direct sum of two line bundles $L',L''$.

Consider the subsheaf $F$ of $E$ consisting of sections of $L'$ on
$C_i$ vanishing at the nodes extended by zero outside $C_i$. Then,
the semistability condition is
$${\chi (L')-s\over w_i}={\chi(F)\over w_i}\le {\chi (E)\over 2}$$
where $s$ denotes the number of nodes in the component.

Similarly, consider the subsheaf $G$ consisting of sections of $L''$
on $C_i$ and sections of $E_{|C_j}$ that glue with $L''$ on the
components other than $C_i$. Then,
$${\chi (G)\over 2-w_i}={\chi (E)-\chi (L')\over 2-w_i}\le
{\chi (E)\over 2}.$$ The two inequalities put together give
$$w_i{\chi(E)\over 2}\le \chi (L')\le w_i{\chi (E)\over 2}+s.$$
The same inequality holds for $L''$. By assumption, $w_i\chi $ is
not an integer, hence, a) follows.

Assume now that the sheaf is not locally free at a node $P$ and that
the torsion-free sheaf is isomorphic at the node to ${\mathcal
M}_P^2$. Denote by $C'_1, C'_2$ the two connected components of
$C-\{ P\} $. Consider the subsheaf $F$ of $E$ consisting of sections
of $E$ on $C'_1$ vanishing at $P$ extended by zero on $C'_2$. Write
$w'_1\ (w'_2)$ for the sum of the $w_i$ corresponding to $C'_1\
(C'_2)$. Then, the semistability condition is
$${\chi (E_{|C'_1})\over 2w'_1}={\chi(F)\over 2w'_1}\le {\chi (E)\over 2}.$$
 Similarly, consider the subsheaf $G$ consisting of sections of
$E$ on $C'_2$ vanishing at $P$ extended by zero on $C'_1$. Then,
$${\chi (E_{|C'_2})\over 2- 2w'_1}={\chi(G)\over 2w'_2}\le {\chi (E)\over 2}.$$
As in this case, $\chi (E)=\chi (E_{|C'_1})+\chi (E_{|C'_2})$, the
two inequalities put together give
$$w'_1\chi(E)\le \chi (E_{C'_1})\le w'_1\chi (E)$$
and this is impossible as, by assumption $w'_1\chi (E)$ is not an
integer. This proves b).
\end{pf}

\end{section}

\begin{section}{Limits on families with singular fiber}

We now want to see how to apply the results above to our situation.
Assume that we have family of curves $\pi :{\mathcal C}\rightarrow
S$ with $S$ one-dimensional such that the fiber over $t_0$ is a
curve of the form \ref{ccgg} while the generic fiber is
non-singular. Note that the structure of the special fiber is
preserved by finite base-change and resolution of singularities.

Assume that we have a family of rank two vector bundles ${\mathcal
E}$ on ${\mathcal C}-{\mathcal C}_0\rightarrow S-\{ t_0\}$ such that
the determinant of the restriction to the fibers is the canonical
line bundle on the fibers. Choose a line-bundle ${\mathcal L}$
defined on ${\mathcal C}\rightarrow S$ such that the
Euler-Poincar\'e characteristic $\chi$ of ${\mathcal E}\otimes
{\mathcal L}$ restricted to each fiber is non-zero. Choose a
polarization $\{w_i\}$ associated to the components $C_i$ of
${\mathcal C}_{t_0}$ so that $\sum _{i\in I} w_i\chi \notin {\bf Z}$
for all $I\subset \{1,\ldots M\}$. Using the arguments of Seshadri
(cf \cite{S}) , replacing a single curve by a family with enough
sections, one has a moduli space of $\{ w_i\}$-semistable
torsion-free sheaves on ${\mathcal C}\rightarrow S$. By properness,
${\mathcal E}$ gives rise to a map from $S$ to this moduli space.
There exists an \'etale covering of the moduli space over which a
universal family exists. Therefore, up to replacing ${\mathcal
C}\rightarrow S$ with a suitable covering, we can assume that
${\mathcal E}$ can be extended to a semistable torsion-free sheaf on
the family of curves that we shall still denote by ${\mathcal E}$.
If the sheaf is not locally free on ${\mathcal C}_{t_0}$, one can
blow up some of the nodes of ${\mathcal C}_{t_0}$ to obtain a new
surface
$$\begin{matrix}{\mathcal C}'&\rightarrow & {\mathcal C}\\
\downarrow & & \downarrow \\
S'& \rightarrow & S\\
\end{matrix}$$
on which the sheaf ${\mathcal E}$ can be extended to a sheaf
${\mathcal E}'$ that is locally free. In this process, one adds a
few rational curves among the nodes of ${\mathcal C}_{t_0}$ and
therefore, the central fiber is still of the type given in
\ref{ccgg}. Moreover, from \ref{stable} b), the restriction of the
vector bundle to these new components is either trivial or of the
form ${\mathcal O}\oplus {\mathcal O}(1)$ (see \cite{G} , section 5
or \cite{X} 1.7). Denote still with ${\mathcal L}$ the pull-back of
the original ${\mathcal L}$ to the new families obtained by
base-change and blow-up. Now ${\mathcal E}'\otimes \pi ^*({\mathcal
L}^{-1})$ restricted to the central fiber can be taken as the limit
of the original ${\mathcal E}$ and satisfies the conditions in
\ref{stable}.

\bigskip

We assume in what follows that we are in the framework we just
described, namely we have a family of curves ${\cal C}\rightarrow S$
where $S$ is the spectrum of a discrete valuation ring. Denote by
$\eta$ the generic point in $S$, by $t$ a generator of the maximal
ideal of the ring and by $\nu$ the discrete valuation.

We assume that ${\cal C}_{\eta}$ is non-singular while the central
fiber is a chain of curves as described in \ref{ccgg}.

 We want to consider linear series for vector bundles of
rank two. On the central fiber, one obtains then a limit linear
series in the sense of \cite{duke}. For the convenience of the
reader, we shall reproduce the definition here.

\begin{Thm}\label{lls} {\bf Definition. Limit linear series}
A limit linear series of rank $r$, degree $d$ and dimension $k$ on
a chain of $M$ (not necessarily rational and elliptic) curves
consists of data I,II below for which data III, IV exist
satisfying conditions a-c.

I) For every component $C_i$, a vector bundle $E_i$ of rank $r$
and degree $d$ and a $k$-dimensional space of sections $V_i$ of
$E_i$.

II) For every node obtained by gluing $Q_i$ and $P_{i+1}$ an
isomorphism of the projectivisation of the fibers $(E_i)_{Q_i}$
and $(E_{i+1})_{P_{i+1}}$

III) A positive integer a

IV) For every node obtained by gluing $Q_i$ and $P_{i+1}$, bases
$s^t_{Q_i}, s^t_{P_{i+1}}, t=1...k$ of the vector spaces $V_i$ and
$V_{i+1}$

Subject to the conditions

a) $\sum _{i=1}^M d_i-r(M-1)a=d$

b) The orders of vanishing at $P_{i+1},Q_i$ of the sections of the
chosen basis satisfy $ord_{P_{i+1}}s_{i+1}^t+ord_{Q_{i}}s_{i}^t\ge
a$

c) Sections of the vector bundles $E_i(-aP_i), E_i(-aQ_i)$ are
completely determined by their value at the nodes.
\end{Thm}

For the reader who is unfamiliar with the concept of limit linear
series, this is to be understood in the following way: On a family
of curves, fix a vector bundle on the generic fiber and extend it
(after some base change and normalizations) to a vector bundle on
the central fiber. This limit vector bundle on the central fiber is
not unique. We could modify it by for example tensoring the bundle
on the whole curve ${\cal C}$   with a line bundle with support on
the central fiber. This would leave the vector bundle on the generic
curve unchanged but would modify the vector bundle on the reducible
curve (by adding to the restriction to each component a linear
combination of the nodes). In this way, for each component $C_i$,
one can choose a version ${\cal E}_i$ of the limit vector bundle
such that the sections of ${\cal E}_i$ that vanish at some of the
nodes of a component $C_j , i\not= j$ are in fact identically zero
on $C_j$. If a vector bundle ${\mathcal E}$ was $\{
w_i\}$-semistable, these modified versions no longer are. But the
restriction to each component preserves the same structure .
Therefore, results as in \ref{stable} are still valid.

As we are allowed to tensor with line bundles supported on the
components of the central fiber, we can assume that $deg({\cal
E}_{i|C_i})=2g-2+\epsilon _i,\ -1\leq \epsilon _i\leq 1, \sum
\epsilon _i=0$. These ${\cal E}_i$ are not independent of each
other. Each ${\cal E}_{i+1}$ can be obtained from ${\cal E}_i$ as
${\cal E}_{i+1}={\cal E}_i(-aF_i)$ where $F_i$ is the union of the
components $C_{i+1},...C_M$. In our case, $a$ can be taken to be
$g-1$

Of special interest to us is the limit of the canonical linear
series . This has rank one and is considered in \cite{EH}
\cite{W}. The restriction of the corresponding line bundle to a
component $C_i$ is given by
$$\omega_{i|C_i}= {\cal O}(2(\sum_{k\le i}g(C_k)-1)P_i+2(g-(\sum_{k\le
i}g(C_k)))Q_i).$$ In the case of a rational curve, this can be
written more simply as ${\cal O}(2g-2)$.

Assume now that we have a limit linear series with canonical
determinant. Write  ${\cal E}_{i|C_i}=E_i$ and write
$degE_i=2g-2+\epsilon _i,\ \  -1\le \epsilon_i\le  1$ as above.

Then,
$$detE_i((\sum_{j<i}\epsilon _j)P_i-(\sum_{j\le
i}\epsilon_j)Q_i)=$$ $$={\cal O}(2(\sum_{j\le
i}g(C_j)-1)P_i+2(g-(\sum_{j\le i}g(C_j))Q_i)$$

Hence,
$$det E_i={\cal O}((2(\sum_{j\le i}g(C_j)-1)-\sum_{j<i}\epsilon_j)P_i+(2(g-(\sum_{j\le
i}g(C_j))+\sum_{j\le i}\epsilon _j))Q_i).$$
 As each $E_i$ is either
semistable or a direct sum of two line bundles whose degrees differ
in one unit, the sections of $E_i(-(g-1)P_i, ),E_i(-(g-1)Q_i)$ are
completely determined by their values at the nodes. We then take
$a=g-1$ in the definition of limit linear series and the condition
$$\sum_idegE_i-2a(M-1)=2g-2$$
is satisfied.

The orders of vanishing at a point $P_i$ (or $Q_i$) of a linear
series will be denoted by $(a_j(P_i))$,(resp $( a_j(Q_i)), \
j=1...k$. Note that each $a_j$ will appear at most twice and when
this happens, there is a two dimensional space of sections with this
vanishing at $P_i$(resp $Q_i$).
\end{section}

\begin{section}{Vanishing at the nodes of elements in the kernel}

 The following results are analogous to 1.2,1.3 of \cite{EH}.
  The proof of the first is almost identical to the one in \cite{EH} and is omitted.

  Here $t$
denotes the parameter in the discrete valuation ring and $\nu $ its
valuation. Choose a component $C_i$ and denote by ${\mathcal E}_i$
as in the previous section, the vector bundle on $\pi :{\mathcal
C}\rightarrow S$ whose restriction to all components of the central
fiber  except $C_i$ has trivial sections.
\begin{Thm}
{\bf Lemma}.  For every component $C_i$, there is a basis $\sigma_j,
j=1...k$ of $\pi_*{\cal E}_i$ such that

a) $ord_{P_i}(\sigma_j)=a_j(P_i)$

b) for suitable integers $\alpha_j$, $t^{\alpha_j}\sigma _j$ are a
basis of $\pi _*({\cal E}_{i+1})$
\end{Thm}
\begin{Thm}\label{vansigma}
{\bf Proposition} Let $\sigma_j$ be a basis of $\pi_*{\cal
E}_{C_i}$such that $t^{\alpha _j}\sigma _j$ is a basis of
$\pi_*({\cal E}_{i+1})$. Then , the orders of vanishing of the
$\sigma_j$ at the nodes satisfy

a) $ord_{P_i}(\sigma_j)\le g-2-ord_{Q_i}\sigma_j\le \alpha _j-1\le
ord_{P_{i+1}}t^{\alpha_j}\sigma_j-1$ if $E_i$ is an indecomposable
vector bundle of degree $2g-3$.

b) $ord_{P_i}(\sigma_j)\le g-1-ord_{Q_i}\sigma_j\le \alpha _j\le
ord_{P_{i+1}}t^{\alpha_j}\sigma_j$ if $E_i$ is an indecomposable
vector bundle of degree $2g-1$, a direct sum of two line bundles of
degree $g-1$, an indecomposable vector bundle of degree $2g-2$ or a
direct sum of two line bundles of degrees $g-1,g-2$.

c) $ord_{P_i}(\sigma_j)\le g-ord_{Q_i}\sigma_j\le \alpha _j+1\le
ord_{P_{i+1}}t^{\alpha_j}\sigma_j+1$ if $E_i$ is a direct sum of two
line bundles of degrees $g-1,\ g$.

Moreover, if equality holds, then $\sigma_j$ vanishes only at
$P_i,Q_i$ as a section of $E_i$.
\end{Thm}

\begin{pf} The proof  is similar to the proof of
the analogous result in Prop 1.1 in \cite{EH}.

 For a line bundle of
degree $d$ , the sum of the vanishing at two given points is at most
$d$ and if equality holds, the section does not vanish at any other
point.

Similarly (see \cite{A}), for an indecomposable vector bundle of
degree $2d+1$ or $2d$ on an elliptic curve, a section vanishes at
two given points with orders adding to at most $d$ and again, if
equality holds, the section does not vanish anywhere else. This
gives the first inequality in each of a,b,c above.

For the second inequality, note that $t^{\alpha _j}\sigma_j$ is a
section of $\pi _*({\mathcal E}_{i+1})=\pi _*({\mathcal
E}_{i}(-(g-1)F_{i+1}))$ where $F_{i+1}$ represents the divisor
$C_{i+1}\cup ..\cup C_M$. Therefore, $\sigma _{j|C_i}$ vanishes as a
section of $\pi _*({\mathcal E}_{i})$ to order at least $g-1-\alpha
_j$ along $F_{i+1}$. Hence, $\sigma _{j|C_i}$ vanishes to order at
least $g-1-\alpha _j$  at $Q_i$. This gives the second of the
inequalities.

For the last inequality, use the fact that ${\mathcal
E}_{i+1}={\mathcal E}_i$ locally along $C-F_{i+1}$. Hence,
$t^{\alpha _j}\sigma _j$ vanishes as a section of $\pi _*({\mathcal
E}_{i+1})$ to order at least $\alpha _j$ on $C-F_{i+1}$. Hence,
$t^{\alpha _j}\sigma _j$ vanishes to order at least $\alpha _j$ at
$P_{i+1}$.

\end{pf}

 \begin{Thm} \label{Remark} {\bf Remark} \end{Thm} The inequalities above can be
 equalities only if there is a section whose order of vanishing at
 $P_i,Q_i$ adds up to the maximum possible. We describe next when
 this happens for each of the possible structures of the restriction
 of the vector bundle to $C_i$. We are assuming that  $C_i$ is
 elliptic and $P_i,\ Q_i$ are generic and   in particular ${\mathcal
 O}(P_i-Q_i)$ is not a torsion point of the Jacobian.

 If $E$ is indecomposable of
degree $2(g-1)+1$ (resp. $2(g-1)-1$), there is a finite number of
sections (up to a constant) such that
$ord_{P_i}(s)+ord_{Q_i}(s)=g-1$ (respectively $g-2$)) . In fact, for
each possible vanishing $a_i$ at $P_i, 0\le a_i\le g-1$ (resp. $0\le
a_i\le g-2$) there is one such section and no two of them are
sections of the same line bundle. This follows from the fact that
$h^0(E_i(-aP-(d-a)Q))=1$ if $deg(E_i)=2d+1$.

In the case of an indecomposable vector bundle of degree $2(g-1)$,
there is at most one section with sum of vanishings at the nodes
being $g-1$.

 If $E=L_g\oplus L_{g-1}$,
there is at most one section that vanishes with maximum sum of
vanishings $g$. There is at most one section of $L_{g-1}$ such that
the sum of vanishing at the nodes is $g-1$. For every value $a$,
there is one section of $L_g$ with vanishing $a$ at $P$ and $g-1-a$
at $Q$.

If $E$ is the direct sum of two line bundles of degree $g-1$, there
are at most two independent sections with maximum vanishing at the
nodes adding up to $g-1$.

\bigskip

We now want to define the order of vanishing of a section $\rho$ of
$S^2(\pi_*({\cal E}_i))$ as follows: Let $\sigma_j$ be a basis of
$\pi_*{\cal E}_i$ such that their orders of vanishing are the orders
of vanishing of the linear series  at $P_i$. Write $\rho =\sum_{j\le
l}f_{jl}(\sigma_j\otimes \sigma_l+\sigma_l\otimes \sigma_j)$ with
$f_{jl}$ functions on the discrete valuation ring.

\begin{Thm}
{\bf Definition} We say $ord_{P_i}\rho \ge \lambda$ if for every
$j,k$ with $f_{j,k}(P_i)\not=0$,
$ord_{P_i}\sigma_j+ord_{P_i}\sigma_k\ge \lambda$
\end{Thm}

 As in
\cite{EH} p.278, one can identify $\pi_*{\cal E}_i$ and $S^2(\pi_*
{\cal E}_i)$ with submodules of $\pi_*{\cal E}_{\eta}$ and
$S^2(\pi_* {\cal E}_{\eta})$ and also with submodules of
$\pi_*{\cal E}_{i+1}$ and $S^2(\pi_* {\cal E}_{i+1})$. Let
$$\rho \in S^2(\pi_* {\cal E}_{\eta}).$$
One can then find a unique value $\beta_i$ such that
$$\rho_i =t^{\beta_i}\rho \in S^2(\pi_* {\cal E}_i)-tS^2(\pi_* {\cal
E}_i).$$

Then for  $\alpha^i=\beta_{i+1}-\beta_i$, one has
$$\rho_{i+1} =t^{\alpha^i}\rho_i \in S^2(\pi_* {\cal E}_{i+1})-tS^2(\pi_* {\cal
E}_{i+1}).$$

The following proposition follows immediately from the definitions:

\begin{Thm}
{\bf Proposition}. \label{ordrho} Fix a component $C_i$. Assume
that $\rho=\sum
 f_{jl}(\sigma_j\otimes\sigma_l+\sigma_l\otimes\sigma_j)$ where
 the $\sigma_j$ are a basis of $\pi_* {\cal E}_i$ such that
 $t^{\alpha_j}\sigma_j$ is a basis of $\pi_*{\cal E}_{i+1}$
$$t^{\alpha^i}\rho \in S^2(\pi_*({\cal E}_i))-tS^2(\pi_*({\cal
E}_i)).$$ Then $$ord_{P_i}(\rho )=min_{\{j,l|
\nu(f_{jl})=0\}}ord_{P_i}(\sigma_j)+ord_{P_i}(\sigma _l)$$
$$\alpha ^i=max_{\{ j,l\} }(\alpha _j+\alpha_l-\nu(f_{jl})).$$
\end{Thm}

We now assume that the kernel of the Petri map is non-trivial on the
generic curve. We can then find a section $\rho _{\eta}$ in the
kernel of the Petri map over the generic point. As above, we can
find a $\rho_i$ for each $i, 1\le i\le M$, in the kernel of the map
$$S^2(\pi_*({\cal E}_i))\rightarrow (\pi_*(S^2{\cal E}_i))$$
with $\rho_i \notin tS^2(\pi_* {\cal E}_i)$.

\begin{Thm}\label{Ld+Ld} {\bf Proposition}
Let $C_i$ be an elliptic curve such that the restriction of the
vector bundle to $C_i$ is the direct sum of two line bundles of
degree $g-1$. Then,
$$ord_{P_{i+1}}(\rho_{i+1})\ge ord_{P_i}(\rho_i)+1$$
Equality above implies that the terms of $\rho_i$ that give the
vanishing at $P_i$ can be written as
$$(\sigma '\otimes \bar \sigma '+ \bar \sigma '\otimes\sigma ')
+ (\sigma ''\otimes \bar \sigma ''+ \bar \sigma ''\otimes\sigma
'')$$ with $\sigma ',\sigma ''$ being  two  independent sections
that vanish at $P_i, Q_i$ with orders adding up to $g-1$ (and in
particular equality implies that these sections exist). Moreover, if
$$(b)\ \ \ ord_{P_i}(\rho)\ge 2(\sum_{j<
i}g(C_j))-\sum_{j<i}\epsilon_j-1=\lambda _i-1$$ then
$$ord_{P_{i+1}}(\rho _{i+1})\ge min(ord_{P_i}(\rho_i )+2,\lambda _i+2) .$$
 When the inequalities for the
order of vanishing at $P_{i+1}$ are equalities, the terms of
$\rho_{i+1}$ that give the minimum vanishing at $P_{i+1}$ glue with
the terms of $\rho _i$ that give the minimum vanishing at $P_i$.
\end{Thm}
\begin{pf}
Assume $C_i$ is an elliptic curve and the vector bundle on $C_i$ is
of the form $L'\oplus L''$, where $L',L''$ are line bundles of
degree $g-1$. Then,
$$L'\otimes L''= {\cal O}((2((\sum_{k\le
i}g(C_k))-1)-\sum_{k<i}\epsilon_k)P_i+ (2(g-(\sum_{k\le
i}g(C_k))+\sum_{k\le i}\epsilon_k)Q_i).$$ For ease of notation, we
write this line bundle as
$${\cal O}(\lambda _i P_i+\mu_i Q_i),
\lambda _i+\mu_i=2(g-1).$$

Let $\rho =\sum f_{j,l}(\sigma_j\otimes \sigma _l+\sigma_l\otimes
\sigma _j)$ be an element in the kernel of the Petri map written in
terms of a basis of $\pi_*{\cal E}_i$ such that the $t^{\alpha
_j}\sigma _j$ is a basis of $\pi_*{\cal E}_{i+1}$. Note that a
single  element of the form $\sigma_j\otimes
\sigma_l+\sigma_l\otimes \sigma_j$ is not in the kernel of the Petri
map. Hence, any $\rho$ in the kernel has in its expression at least
two summands of this type with the $\sigma$ 's appearing there being
independent (we could have something like $c_{11}\sigma_1\otimes
\sigma_1 +c_{22}\sigma_2\otimes \sigma_2$ but not something like
$c_{12}(\sigma_1\otimes \sigma_2 +\sigma_2\otimes
\sigma_1)+c_{13}(\sigma_1\otimes \sigma_3 +\sigma_3\otimes
\sigma_1)$ as the latter could be written as $\sigma_1\otimes
(c_{12}\sigma_2
+c_{13}\sigma_3)+(c_{12}\sigma_2+c_{13}\sigma_3)\otimes \sigma_1$.

Note that $ord_{P_i}(\sigma_j)\le \alpha_j\le ord_{P_{i+1}}t^{\alpha
_j}\sigma_j$ with the first inequality being an equality for at most
two independent sections $\sigma ' , \sigma ''$ that vanish at
$P_i,Q_i$ with orders adding up to $g-1$ (see \ref{vansigma}
\ref{Remark}).

Hence, $ord_{P_{i+1}}(\rho _{i+1}) \ge ord_{P_i}(\rho _i)+1$ and
$ord_{P_{i+1}}(\rho _{i+1}) \ge ord_{P_i}(\rho _i)+2$ except in the
case when the  expression of $\rho$ that gives the vanishing at
$P_i$ contains exactly two summands and one of the sections
appearing in each summand is one of the two sections $\sigma ',\
\sigma ''$ named above. So, we can write the restriction of $\rho$
to $C_i$ as
$$\rho =(\sigma '\otimes \bar \sigma '+ \bar \sigma '\otimes\sigma ')+
(\sigma''\otimes \bar \sigma''+ \bar \sigma ''\otimes\sigma
'')+...$$ where the dots stand for sections with higher vanishing at
$P_i$. Using the decomposition of $E$ into a direct sum and up to
replacing $\bar \sigma '$ and $\bar \sigma ''$ by linear
combinations of themselves, assume that $\sigma '=(\sigma _1',0)$
and $\sigma ''=(0,\sigma _2'')$ . Write then $\bar \sigma
'=(\bar\sigma_1', \bar \sigma_2')$ (and similarly for $\bar \sigma
''$) in terms of the decomposition of $E_i$. Using that $\rho$ is in
the kernel of the Petri map, we obtain $\bar\sigma'_1=0,\ \bar
\sigma''_2=0$. Hence, the order of vanishing of $\bar \sigma '$ is
in fact the order of vanishing of $\bar \sigma'_2$. From the fact
that $\sigma _1'\sigma _2''$ is a section of $L'L"$ vanishing to
order adding up to $2g-2$ between the two nodes,
 $ord_{P_i}(\sigma ' \sigma '')=\lambda _i$ and
$ord_{P_i}(\sigma'\bar \sigma ')\not=ord_{P_i}(\sigma '\sigma '').$
Assume now that $ord_{P_i}(\rho_i)\ge \lambda_i$. Then
$$\lambda _i
\le ord_{P_i}(\rho_i)=ord_{P_i}(\sigma')+ord_{P_i}(\bar\sigma
')\not=
$$
$$\not= ord_{P_i}(\sigma')+ord_{P_i}(\sigma ")=\lambda_i.$$ Hence
$$\lambda _i+1\le
ord_{P_i}(\rho _i)=ord_{P_i}(\sigma ')+ord_{P_i}(\bar \sigma ').$$

On the other hand by \ref{vansigma} \ref{ordrho},
$$ord_{P_{i+1}}(t^{\alpha ^i}\rho _i )\ge ord_{P_{i+1}}(t^{\alpha '}\sigma
' )+ord_{P_{i+1}}(t^{\bar\alpha '}\bar \sigma ')\ge ord_{P_i}(\sigma
') +ord_{P_i}(\bar \sigma ')+1.$$ The two inequalities together
conclude the proof of the result when $ord_{P_i}(\rho_i)\ge \lambda
_i$.

Assume now that $ord_{P_i}(\rho_i)= \lambda _i-1$ and we want to
show that also in this case $ord_{P_{i+1}}(\rho_{i+1})\ge
ord_{P_i}(\rho_i) +2$.

We can again assume that the terms in  $\rho$ that gives the
minimum vanishing at $P_i$ are of the form
$$(\sigma '\otimes \bar \sigma '+ \bar \sigma '\otimes\sigma ')
+ (\sigma''\otimes \bar \sigma ''+ \bar \sigma ''\otimes\sigma '')$$

 and that
$\alpha^i=\alpha'+\bar \alpha '=\bar \alpha ''+\alpha ''
=ord_{P_i}(\rho_i) +1$ otherwise using \ref{vansigma} and
\ref{ordrho}, we would be done. Note now that $h^0(C_i,L'\otimes
L''(-(\lambda _i-1)P_i-\mu_iQ_i))=1$ because the degree of this line
bundle is one. Also $h^0(C_i,L'\otimes L''(-\lambda
_iP_i-\mu_iQ_i))=h^0(C,{\cal O})=1$ while $h^0(L'\otimes
L''(-(\lambda _i-1)P_i-(\mu_i-1)Q_i))=2$. As
$$H^0(C,L'\otimes L''(-\lambda _iP_i-\mu_iQ_i))\subset H^0(C,L'\otimes
L''(-(\lambda _i-1)P_i-\mu_iQ_i))$$ these two vector spaces having
the same dimension, coincide and a section of $L'\otimes L''$ that
vanishes at $P_i$ with multiplicity precisely $\lambda_i-1$
vanishes at $Q_i$ with multiplicity at most
$\mu_i-1=2g-2-\lambda_i-1=2g-2-(\lambda_i+1)$. Then a section
gluing with it vanishes at $P_{i+1}$ with order at least
$\lambda_i+1=(\lambda_i-1)+2$. In particular,
$(t^{\alpha'}\sigma')(t^{\bar \alpha '}\bar \sigma')$ and
$(t^{\alpha "}\sigma ")(t^{ \bar \alpha "}\bar \sigma")$ satisfy
the condition. Then,
$$\rho _{i+1}=t^{\alpha^i} \rho_i=(t^{\alpha '}\sigma '\otimes
 t^{\bar \alpha '}\bar \sigma '+
t^{\bar \alpha '}\bar \sigma'\otimes t^{\alpha '}\sigma ')+ (t^{
\alpha "}\sigma "\otimes t^{\bar \alpha "}\bar \sigma "+ t^{\bar
\alpha "}\sigma "\otimes t^{\alpha "}\bar \sigma ")$$ vanishes at
$P_{i+1}$ with the required multiplicity.

If there is equality, $\alpha ^i=\alpha '+\bar \alpha '=\alpha
"+\bar \alpha "$ and the terms in $\rho_{i+1}$ that give the
vanishing at $P_{i+1}$ glue with the terms in $\rho_i$ that give the
vanishing at $P_i$.
\end{pf}

\begin{Thm}
{\bf Proposition} Assume that $ E_i$ is an indecomposable vector
bundle of degree $2(g-1)$. Then
$$ord_{P_{i+1}}(\rho _{i+1})\ge ord_{P_i}(\rho _i)+2.$$
If the inequality is an equality, the terms in $\rho_{i+1}$ that
give minimum vanishing at $P_{i+1}$ glue with the terms in
$\rho_i$ that give minimum vanishing at $P_i$.
\end{Thm}

\begin{pf} In this case there is at most one section $\sigma$
that vanishes at $P_i$ and $Q_i$ with multiplicity adding up to
$g-1$ (see  \ref{Remark}). As an element of the form $\sigma \otimes
\bar \sigma + \bar \sigma \otimes\sigma $ is not in the kernel of
the Petri map,
 in the terms of $\rho_i$ that give minimum vanishing
at $P_i$, $\rho_i =\sum_{jl}(\sigma_j\otimes
\sigma_l+\sigma_l\otimes\sigma_j)+...$ for at least one pair $\{
\sigma_j, \sigma_l\}$ does not contain $\sigma $. Hence the
inequality follows. If some of the terms on $\rho_{i+1}$ that give
the vanishing at $P_{i+1}$ come from terms on $\rho_i$ with
vanishing at $P_i$ at least $ord_{P_i}\rho_i+1$, then by the same
argument, $ord_{P_{i+1}}(\rho_{i+1})\ge
(ord_{P_{i}}(\rho_{i})+1)+2$. This proves the second statement.

\end{pf}

\begin{Thm}\label{2d+1}
{\bf Proposition} Assume that $ E_i$ is an indecomposable vector
bundle of degree $2(g-1)+\epsilon, \ \epsilon=1,-1$. Let $\rho_i$
be as before Then either
$$(a)\ \ \ ord_{P_{i+1}}(\rho _{i+1} )\ge ord_{P_i}(\rho _i)-\epsilon +2$$
or b) $$ord_{P_{i+1}}(\rho _{i+1}) \ge ord_{P_i}(\rho _i)-\epsilon
+1$$ and the terms in $\rho_i$ that give the vanishing at $P_i$ can
be written in the form
$$\rho=c_{12}(\sigma_1\otimes \sigma _2+\sigma_2\otimes \sigma _1)+
c_{34}(\sigma_3\otimes \sigma _4+\sigma_4\otimes \sigma _3)+...$$
where all the $\sigma_j$ are sections such that the sums of
vanishing at $P_i,Q_i$ is $g-1$ if $\epsilon =1$, $g-2$ if $\epsilon
=-1$ and $\rho$ does not admit an expression with fewer than three
summands $\sigma_j\otimes \sigma_l+\sigma_l\otimes \sigma _j$.

Also, when equality occurs, the terms in $\rho_{i+1}$ that give
minimum vanishing at $P_{i+1}$ glue with the  terms in $\rho_i$
with minimum vanishing at $P_i$.
\end{Thm}
\begin{pf} We prove it in the case $\epsilon =1$, the other being
analogous.

 If $\rho _i=\sum f_{j,l}(\sigma_j \otimes
\sigma_l+\sigma_l\otimes \sigma_j)$ is an element in the kernel,
then there are at least two pairs $j,l$ that give the minimum
vanishing at each point in the above expression. If in one of these
pairs at most one of the sections has  vanishings at $P_i,Q_i$
adding up to $g-1$, using \ref{vansigma},  \ref{ordrho}  above, we
are done. Assume then that in the above expression all the sections
that appear have sum of vanishings at $P_i,Q_i$ equal to $g-1$. We
want to check that then, the expression for $\rho $ has at least
three summands.

Choose a line bundle $L'$ of degree $g-1$ not linearly equivalent to
$aP+(g-1-a)Q, 0\le a\le g-1$. Let $L''$ be the line bundle of degree
$g$ such that $L'\otimes L''=det E$. Then, $E$ can be deformed to a
direct sum $L'\oplus L''$. On $L'\oplus L''$ there are exactly $g$
sections $s _j$ with sum of  vanishing at $P_i,Q_i$ being $g-1$. In
fact $ord_{P_i}(s_j)=j, ord_{Q_i}(s_j)=g-1-j$. As $\deg
(L''(-jP_i-(g-1-j)Q_i))=1$, each one of these sections vanishes at
an additional point $R_j$. The $R_j$ are all different: if $R_j=R_k,
\ j>k$, then from the fact that $s_j, s_k$ are sections of the same
line bundle $L"$, we obtain $(j-k)P_i\equiv (j-k)Q_i$ contradicting
the genericity of the pair $P_i,Q_i$. Assume then
$$c_{i_1i_2}(\sigma_{i_1}\otimes \sigma _{i_2}+\sigma_{i_2}\otimes
\sigma _{i_1})+c_{i_3i_4}(\sigma_{i_3}\otimes \sigma
_{i_4}+\sigma_{i_4}\otimes \sigma _{i_3})$$ is in the kernel of the
Petri map. Then,
$2c_{i_1i_2}s_{i_1}s_{i_2}=-2c_{i_3i_4}s_{i_3}s_{i_4}$. Hence these
two sections vanish at the same points and we obtain $\{
R_{i_1},R_{i_2}\} =\{ R_{i_3},R_{i_4}\} $ which is a contradiction.
\end{pf}
\begin{Thm}
{\bf Proposition} Assume that $C_i$ is elliptic and $ E_i$ is a
direct sum of two line bundles of degree $g-1$ and $g-1+\epsilon ,
\epsilon=1,-1$ respectively . Then either
$$ord_{P_{i+1}}(\rho_{i+1}) \ge ord_{P_i}(\rho_i)-\epsilon +2$$
or $$ord_{P_{i+1}}(\rho_{i+1}) \ge ord_{P_i}(\rho_i)-\epsilon +1$$
and $\rho_i$ can be written in the form
$$\rho_i=\sum f_{jl}(\sigma_j\otimes \sigma _l+\sigma_l\otimes \sigma _j)+
f'(\sigma '\otimes \sigma +\sigma\otimes \sigma ')+f''(\sigma
''\otimes \sigma _i+\sigma_i\otimes \sigma ")$$ where the $\sigma_j,
\sigma_l, \sigma _i,\ \sigma '$ are all sections of the line bundle
of higher degree.

Moreover, if
$$ord_{P_i}(\rho)\ge 2((\sum_{j\le
i}g(C_j))-1)-\sum_{j<i}\epsilon_j-1,$$ then
$$ord_{P_{i+1}}(\rho _{i+1})\ge ord_{P_i}(\rho _i)+2$$
or $\rho$ can be written as before with $f'=f"=0$. When the
inequalities for the order of vanishing at $P_{i+1}$ are equalities,
the terms in $\rho_{i+1}$ that give the minimum vanishing at
$P_{i+1}$ glue with the terms in  $\rho _i$ that give the minimum
vanishing at $P_i$.

\end{Thm}
\begin{pf} We assume $\epsilon=1$, the case $\epsilon =-1$ being
similar.

Write $E_i=L'\oplus L"$ where $L'$ is a line bundle of degree $g$
and $L"$ is a line bundle of degree $g-1$. There are then at most
one section of $L'$ and one section of $L"$ that vanish only at
$P_i,Q_i$ with sum of vanishing at these two points being $g,g-1$
respectively. For these two sections, one has
$$ord_{P_i}(\sigma ')\le g-ord_{Q_i}(\sigma ')\le \alpha'+1\le
ord_{P_{i+1}}(t^{\alpha '}\sigma ')+1$$
$$ord_{P_i}(\sigma '')\le g-1-ord_{Q_i}(\sigma '')\le \alpha ''\le
ord_{P_{i+1}}t^{\alpha ''}\sigma ''.$$
 For any other section of $L'$ one has
$$ord_{P_i}(\sigma _j')\le g-1-ord_{Q_i}(\sigma _j')\le \alpha'_j \le
ord_{P_{i+1}}t^{\alpha '_j}\sigma _j'$$

while for other sections $\sigma_j ''$ one has
 $$ord_{P_i}(\sigma _j'')\le g-2-ord_{Q_i}(\sigma _j'')\le \alpha _j''-1\le
ord_{P_{i+1}}t^{\alpha _j''}\sigma _j''-1.$$

Write now
$$\rho_i=\sum f_{jl}(\sigma_j\otimes \sigma _l+\sigma_l\otimes \sigma
_j).$$

As $\rho _i$ is in the kernel of the Petri map, not all terms in
this expression can contain $\sigma '$. Hence, $ord_{P_{i+1}}(\rho
_{i+1})\ge ord_{P_i}(\rho _i)$.

If this expression has some $\sigma_j$ or $\sigma_l$ that are not
paired with  linear combinations of $\sigma', \sigma''$ or the
sections $\sigma'_j$ , then the vanishing of $\rho _{i+1}$ at
$P_{i+1}$ increases by at least one unit with respect to that at
$P_i$. If this increase does not take place, using the decomposition
of $E_i$ into direct sum, we can write
$$\rho _i=f'((\sigma',0)\otimes (\sigma' _1, \sigma"_1)+(\sigma' _1, \sigma"_1)\otimes (\sigma
',0))+$$
$$f"((0,\sigma '')\otimes (\sigma' _2,\sigma''_2)+(\sigma '_2,\sigma''_2)\otimes (0,\sigma
''))+$$
$$+\sum f_{jl}((\sigma_j,0)\otimes (\sigma _l,0)+(\sigma_l,0)\otimes (\sigma
_j,0).$$

 From the condition of this section being in the kernel, one
 obtains $$\sigma''_2=0,\  f'\sigma ' \sigma ''_1+f''\sigma'' \sigma'_2 =0.$$

If $f', \ f''= 0$, we
 are in the case when $\rho$ can be written only in terms of the
 sections of the line bundle of higher degree.
If $f',\ f''$ are not both zero, then both must be non-zero by the
equation above and we proceed as in the case of the sum of two line
bundles of the same degree to get the result.

\end{pf}
\begin{Thm}\label{racional}
{\bf Proposition} If $C_i$ is a rational curve and $E_i$ is the
direct sum of two line bundles of degree $g-1$, then
$$ord_{P_{i+1}}(\rho_{i+1})\ge ord_{P_i}(\rho_i).$$
If $C_i$ is a rational curve and $E_i$ is a direct sum of two line
bundles of degree $g-1, g-1+\epsilon$, then, either
$$ord_{P_{i+1}}(\rho_{i+1})\ge ord_{P_i}(\rho_i)-\epsilon$$ or
$$ord_{P_{i+1}}(\rho_{i+1})\ge ord_{P_i}(\rho_i)-\epsilon -1$$ and $\rho_i$
can be written in the form
$$\rho_i
=\sum f_{jl}(\sigma_j\otimes \sigma_l+\sigma_l\otimes \sigma_j)$$
where all of the $\sigma _j, \sigma _l$ are sections of the line
bundle of higher degree.
\end{Thm}
The proof of this last claim is similar to the previous ones and
is left to the reader.

\begin{Thm}\label{seccionsdep} {\bf Proposition} Let $C_i$ be elliptic and the
restriction of $E$ to $C_i$ a direct sum of two line bundles of
degree $g-1$. If in the expression of the restriction of $\rho$ to
$C_i$ all the sections of $E$ that appear are linearly dependent at
$P_i$, then
$$ord_{P_{i+1}}(\rho_{i+1})\ge ord_{P_i}(\rho_i)+2$$
and equality above implies that in the expression of $\rho_i$
there are at least three terms or that the expression for $\rho
_i$ can be written in terms of sections of a line subbundle of
$E_i$.
\end{Thm}
\begin{pf}
The inequality has already been proved in \ref{Ld+Ld}.

Assume now that there is equality and the expression for $\rho _i$
has precisely two terms. Because of the condition on the fibers at
$P_i$ at most one of the sections appearing in the expression for
$\rho _i$ has vanishing at the nodes adding up to $g-1$. Call this
section $\sigma '$ and  the line bundle it generates $L$. Write
$$\rho _i=(\sigma '\otimes \bar \sigma '+\bar \sigma '\otimes \sigma ')+
(\sigma ''\otimes \bar \sigma ''+ \bar \sigma ''\otimes \sigma
'').$$ Then $\sigma '',\ \bar \sigma ''$ have order of vanishing at
$P_i,Q_i$ adding up to $g-2$. Therefore they are completely
determined by their direction at $P$. This implies that $\sigma '',\
\bar \sigma ''$ are also sections of the line bundle $L$ and hence
so is $\bar \sigma '$ by the condition of $\rho$ being in the kernel
of the Petri map. So the result is satisfied in this case.

Assume now that none of the sections involved vanishes to the
maximum order $g-1$ between the two nodes. Write the sections
appearing in the expression in terms of the decomposition of $E_i$
as
$$((\sigma'_ 1,\sigma '_2)\otimes (\bar \sigma'_ 1, \bar \sigma '_2)+
(\bar \sigma'_ 1,\bar \sigma '_2)\otimes (\sigma'_ 1,\sigma '_2))+$$
$$ ((\sigma''_ 1,\sigma ''_2)\otimes (\bar \sigma''_ 1, \bar \sigma
"_2)+ (\bar \sigma''_ 1,\bar \sigma''_2)\otimes (\sigma''_ 1, \sigma
"_2)) ).$$

Each of the $\sigma $ that appears in the expression above is either
zero or is the unique section of $L$ (resp $L'$) whose order of
vanishing at $P_i$ is a given $a$ and at $Q_i$ is at least $g-2-a$.
Then, the condition of $\rho $ being in the kernel implies that

$$(*)\sigma'_1\bar \sigma'_1+\sigma''_1\bar \sigma''_1=0,\
\sigma'_2\bar \sigma'_2+\sigma''_2\bar \sigma''_2=0, $$
$$\sigma'_1\bar \sigma'_2+\sigma'_2\bar \sigma'_1+\sigma''_1\bar
\sigma''_2+\sigma''_2\bar \sigma''_1=0.$$

 An argument as in the
proof of  \ref{2d+1} implies that
$$\sigma'_1=\lambda \sigma''_1,\
\bar \sigma''_1=-\lambda \bar \sigma'_1\ {\rom or} \
\sigma'_1=\lambda \bar \sigma''_1,\ \sigma''_1=-\lambda \bar
\sigma'_1.$$ Similarly $$ \sigma '_2=\mu \sigma ''_2,\ \bar \sigma
''_2=-\mu \bar \sigma '_2 \ {\rom or}  \   \sigma '_2=\mu \bar
\sigma ''_2,\ \sigma ''_2=-\mu \bar \sigma '_2 .$$ Let us assume
that we are in the first case for both equations (the remaining
cases are treated similarly). Substituting in the third equation
above, we obtain
$$(\lambda -\mu)(\sigma''_1\bar \sigma'_2-\sigma''_2\bar
\sigma'_1)=0.$$ If $\lambda =\mu$, the original section was trivial
against the assumption. If $\lambda \not= \mu$, then
$(\sigma''_1\bar \sigma'_2-\sigma''_2\bar \sigma'_1)=0$ and an
argument as in the proof of \ref{2d+1} implies that
$$\sigma''_1=\nu \sigma''_2,\
\bar \sigma'_1=\nu \bar \sigma'_2\ {\rom or} \ \sigma''_1=\nu \bar
\sigma'_1,\ \sigma''_2=\nu \bar \sigma'_2.$$

In both cases, the original section was trivial against the
assumptions.
\end{pf}

\begin{Thm}{\bf Proposition} Let $C$ be an elliptic curve and assume that $E_i$ is
an indecomposable vector bundle of degree $2g-1$. Consider the
sections of the vector bundle that vanish at $P$ with multiplicity
$a$ and at $Q$ with multiplicity $g-1-a$. Then, the fibers of these
sections at $P$ are different except maybe for some pairs $a_1,a_2$
and if this happens then $det E_i={\mathcal
O}((a_1+a_2)P+(2g-2-a_1-a_2)Q)$.
\end{Thm}

\begin{pf} Consider two line subbundles  of $E_i$ generated by these
sections, say
 $${\mathcal O}(a_1P+(g-1-a_1)Q)\rightarrow E_i$$
 $${\mathcal O}(a_2P+(g-1-a_2)Q)\rightarrow E_i.$$
 Assume $a_1>a_2$ and write $a'=a_1-a_2,\ b'=g-1-a_2,\ E'=E_i(-a_2P)$ tensoring
 the spaces in the two maps above with ${\mathcal O}(-a_2P)$, we
 obtain
 $$\begin{matrix} 0&\rightarrow &{\mathcal
 O}(a'P+(b'-a')Q)&\rightarrow& E'&\rightarrow & L&\rightarrow &0\\
 & & & & \uparrow & & & & \\
  & & & & {\mathcal O}(b'Q) & & & & \\
  \end{matrix}$$
  As the vertical arrow does not factor through ${\mathcal
  O}(a'P+(b'-a')Q)$, it gives rise to a non-zero morphism to $L$
  (which is a line bundle of degree $b'+1$). If this map is zero at
  $P_i$, then $L={\mathcal O}(P_i+b'Q_i)$. In particular $E_i$ is as claimed
  and $a_2$ is determined by $E_i$ and $a_1$.
  \end{pf}

\end{section}
\bigskip

\begin{section}{Proof of the Theorem}
\begin{Thm}{\bf Proposition}
For a section $\rho$ in the kernel of the Petri map, one has one
of the three options below

a) $$ord _{P_i}(\rho)\ge 2(\sum_{k<i} g(C_k))-\sum_{k<i}\epsilon_k$$

b)
$$ord _{P_i}(\rho)\ge 2(\sum_{k<i}
g(C_k))-\sum_{k<i}\epsilon_k-1$$ and on $C_{i-1}$ the equation for
the terms of $\rho_t$ giving the maximum vanishing at $Q_{i-1}$ is
of the form
$$\sum_{jl}f_{jl}(\sigma_j\otimes \sigma_l+\sigma_l\otimes
\sigma_j)$$ cannot be written with fewer than three summands and
the $\sigma$ appearing in the expression  vanish only at the two
nodes.

c)$$ord _{P_i}(\rho)\ge 2(\sum_{k<i}
g(C_k))-\sum_{k<i}\epsilon_k-1$$ and  on $C_{i-1}$ the equation for
the terms of $\rho_t$ giving the maximum vanishing at $Q_{i-1}$ can
be written in terms of sections of a line subbundle of
$C_{i-k}\cup...\cup C_{i-1}$ of degree ${\deg (E_{i-k})+...+\deg
(E_{i-1})\over 2}+1$

\end{Thm}
\begin{pf} By induction on $i$, the case $i=1$ being clear as
$$2(\sum_{k<1} g(C_k))-\sum_{k<1}\epsilon_k=0$$
Assume the result correct for $i-1$ and show it for $i$. If
condition a) holds, then the result can be obtained from the
propositions in the previous paragraph.

If condition b) holds, and $C_{i}$ is rational or elliptic with
$E$ the direct sum of two line bundles,  then the result follows
also from the previous propositions.

Assume that $C_{i}$ is elliptic and $E_{i}$ is an indecomposable
vector bundle of odd degree. If the vanishing at $P_i$ increases by
at least one unit, we are back to conditions a) or b). Otherwise,
the restriction of $\rho$ to $C_i$ has at least three terms.  With a
suitable choice of gluing of $Q_i,P_{i+1}$ up to two arbitrary
directions can be made to agree. By the independent genericity of
$C_{i-1}, C_i$ not all directions involved in $\rho_i$ will glue
with those in $\rho_{i-1}$. Hence, this is impossible.

Assume now that condition c) holds. If $C_i$ is elliptic and $E_i$
is an indecomposable vector bundle, we reason as in case b). If
$C_i$ is elliptic or rational and $E_i$ is the direct sum of two
line bundles of the same degree, then \ref{Ld+Ld}, \ref{racional}
and \ref{seccionsdep}, above show that either a) or c) holds again.
If $C_i$ is elliptic or rational and $E_i$ is the direct sum of two
line bundles of different degrees, then by stability, the line
bundle of higher degree cannot glue with a line bundle of higher
degree in the previous component. Hence, the result is satisfied
again.

\end{pf}
The Theorem now follows from the statement above: on the last
component (that we can assume to be rational), the vector bundle is
the direct sum of two line bundles of the same degree $g-1$ and
$\sum \epsilon_j=0$. The proposition above implies that $\rho$
vanishes at $P_g$ with multiplicitly at least $2g-1$ and this is
impossible as this multiplicity is at most $2g-2$. Hence $\rho$
cannot exist and the Petri map is injective.

\end{section}


\begin{thebibliography}{cccc}

\bibitem[A] {A} M.Atiyah, {\it Vector bundles over an elliptic
curve} Proc. London Math.Soc.{\bf (3),7} (1957, 414-452.

\bibitem[BF]{BF} A.Bertram, B.Feinberg {On stable rank two bundles
with canonical determinant and many sections}, in "Algebraic
Geometry", ed P.Newstead, Marcel Dekker 1998, 259-269.


\bibitem[EH]{EH} D.Eisenbud, J.Harris, {\it A simpler proof of the
Gieseker-Petri Theorem on special divisors}, Invent.Math. {\bf 74}
(1983), 269-280.

\bibitem[FP]{FP} G.Farkas, M.Popa, {\it Effective divisors on
$\overline{\mathcal M}_g$, curves on $K3$ surfaces, and the slope
conjecture}  J. Algebraic Geom.  {\bf 14}  (2005),  no. 2,
241--267

\bibitem[G]{G} D.Gieseker {\it A degeneration of the moduli space
of vector bundles} J.Diff Geom. {\bf 19}(1984), 173-206.

\bibitem[M1]{M1} S.Mukai, {\it Curves and Brill-Noether Theory},
MSRI publications {\bf 28}, Cambridge University Press 1995,
145-158.

\bibitem[M2]{M2} S.Mukai {\it Non-abelian Brill-Noether Theory and
Fano threefolds}.  Sugaku Expositions  {\bf 14}  (2001),  no. 2,
125--153.

 \bibitem[NS]{NS} D. Nagaraj, C.S.Seshadri, {\it Degenerations of the moduli spaces
of vector bundles on curves. I.}  Proc. Indian Acad. Sci. Math.
Sci.  {\bf 107} (1997),  no. 2, 101--137.

\bibitem[OPP]{OPP} W.Oxbury, C.Pauly, E.Previatto, {\it
Subvarieties of ${\cal SU}_C(2)$ and $2\theta$ divisors in the
Jacobian}, Trans. of the Amer. Math. Soc. {\bf 350
N9},(1998)3587-3617.

\bibitem[P] {P} S.Park {\it Non-emptiness of Brill-Noether loci in
$M(2,K)$} Preprint.

\bibitem[S] {S} C.S.Seshadri {\it Fibr\'es vectoriels sur les courbes
alg\'ebriques} Ast\'erisque 96, Soci\'et\'e Math\'ematique de France
1982.


\bibitem[T1]{duke} M.Teixidor, {\it Brill-Noether Theory for stable vector
bundles} Duke Math J{\bf 62 N2}(1991), 385-400.

\bibitem[T2]{arbre} M.Teixidor, {\it Moduli spaces of semistable vector bundles on
tree-like curves}, Math Ann. {\bf 290} (1991), 341-348.

\bibitem[T3]{canonic} M.Teixidor, {\it Rank two vector bundles with canonical
determinant},   Math. Nachr.  {\bf 265},   (2004), 100--106.

\bibitem[W]{W} G.Welters, {\it A Theorem of Gieseker-Petri type for Prym
varieties}, Ann.scient.E.Norm.Sup. {\bf 4serie, t.18}, 1985,
671-683.

\bibitem[X]{X} H.Xia, {\it Degenerations of moduli of stable
bundles over algebraic curves}, Comp. Math. {\bf 98, n3} (1995),
305-330.

\end{thebibliography}
\end{document}